# A GENERAL THEOREM FOR THE CHARACTERIZATION OF N PRIME NUMBERS SIMULTANEOUSLY

Florentin Smarandache

University of New Mexico

200 College Road

Gallup, NM 87301, USA

1. **ABSTRACT**. This article presents a necessary and sufficient theorem for N numbers, coprime two by two, to be prime simultaneously.

It generalizes V. Popa's theorem [3], as well as I. Cucurezeanu's theorem ([1], p. 165), Clement's theorem, S. Patrizio's theorems [2], etc.

Particularly, this General Theorem offers different characterizations for twin primes, for quadruple primes, etc.

2. **INTRODUCTION**. It is evidently the following:

<u>Lemma 1</u>. Let A,B be nonzero integers. Then:

$AB \equiv 0 \pmod{pB} \Leftrightarrow A \equiv 0 \pmod{p} \Leftrightarrow A/p$ is an integer.

Lemma 2. Let $(p,q) \sim 1$, $(a,p) \sim 1$, $(b,q) \sim 1$. Then:

$A \equiv 0 \pmod{p}$ and $B \equiv 0 \pmod{q} \Leftrightarrow aAq + bBp \equiv 0 \pmod{pq} \Leftrightarrow aA + bBp/q \equiv 0 \pmod{p} \Leftrightarrow aA/p + bB/q$ is an integer.

Proof:

The first equivalence:

We have $A = K_1 p$ and $B = K_2 q$, with $K_1, K_2 \in Z$, hence

$aAq + bBp = (aK_1 + bK_2)pq$.

Reciprocal: $aAq + bBp = Kpq$, with $K \in Z$, it results that $aAq \equiv 0 \pmod{p}$ and $bBp \equiv 0 \pmod{q}$, but from our assumption we find $A \equiv 0 \pmod{p}$ and $B \equiv 0 \pmod{q}$.

The second and third equivalence result from lemma 1.

By induction we extend this lemma to

LEMMA 3: Let $p_1, \ldots, p_n$ be coprime integers two by two, and let $a_1, \ldots, a_n$ be integer numbers such that $(a_i, p_i) \sim 1$ for all i. Then:

$A_1 \equiv 0 \pmod{p_1}, \ldots, A_n \equiv 0 \pmod{p_n} \Leftrightarrow$

$\Leftrightarrow \sum_{i=1}^{n} a_i A_i \cdot \prod_{j \neq i} p_j \equiv 0 \pmod{p_1 \ldots p_n} \Leftrightarrow$

$\Leftrightarrow (P/D) \cdot \sum_{i=1}^{n} (a_i A_i / p_i) \equiv 0 \pmod{P/D}$,

where $P = p_1 \ldots p_n$ and D is a divisor of p, $\Leftrightarrow$

$$\Leftrightarrow \sum_{i=1}^{n} a_i A_i / p_i \text{ is an integer.}$$

3. From this last lemma we can immediately find a

**GENERAL THEOREM**:

Let $P_{ij}$, $1 \leq i \leq n$, $1 \leq j \leq m_i$, be coprime integers two by two, and let $r_1, \ldots, r_n, a_1, \ldots, a_n$ be integer numbers such that $a_i$ be coprime with $r_i$ for all i.

The following conditions are considered:

(i)

$p_{i_1}, \ldots, p_{in_1}$, are simultaneously prime if and only if $c_i \equiv 0 \pmod{r_i}$, for all i.

Then:

The numbers $p_{ij}$, $1 \leq i \leq n$, $1 \leq j \leq m_i$, are simultaneously prime if and only if

(*)  $(R/D) \cdot \sum_{i=1}^{n} (a_i c_i / r_i) \equiv 0 \pmod{R/D}$,

where $P = \prod_{i=1}^{n} r_i$ and D is a divisor of R.

<u>Remark.</u>

Often in the conditions (i) the module $r_i$ is equal to $\prod_{j=1}^{m_i} p_{ij}$, or to a divisor of it, and in this case the relation of the General Theorem becomes:

$$(P/D) \cdot \sum_{i=1}^{n} (a_i c_i / \prod_{j=1}^{m_i} p_{ij}) \equiv 0 \pmod{P/D},$$

where

$$P = \prod_{i,j=1}^{n,m_i} p_{ij} \text{ and } D \text{ is a divisor of } P.$$

Corollaries.

We easily obtain that our last relation is equivalent to:

$$\sum_{i=1}^{n} a_i c_i \cdot (P / \prod_{j=1}^{m_i} p_{ij}) \equiv 0 \pmod{P},$$

and

$$\sum_{i=1}^{n} (a_i c_i / \prod_{j=1}^{m_i} p_{ij}) \text{ is an integer,}$$

etc.

The imposed restrictions for the numbers $p_{ij}$ from the General Theorem are very wide, because if there were two non-coprime distinct numbers, then at least one from these would not be prime, hence the $m_1 + ... + m_n$ numbers might not be prime.

The General Theorem has many variants in accordance with the assigned values for the parameters $a_1, \ldots, a_n$, and $r_1, \ldots, r_m$, the parameter D, as well as in accordance with the congruences $c_1, \ldots, c_n$ which characterize either a prime number or many other prime numbers simultaneously. We can start from the theorems (conditions $c_i$) which characterize a single prime number [see Wilson, Leibniz, Smarandache [4], or Simionov (p is prime if and only if $(p-k)!(k-1)!-(-1)^k \equiv 0 \pmod{p}$, when $p \geq k \geq 1$; here, it is preferable to take $k = \lfloor 1(p+1)/2m \rfloor$, where $\lfloor x \rfloor$ represents the greatest integer number $\leq x$, in order that the number $(p-k)!(k-1)!$ be the smallest possible] for obtaining, by means of the General Theorem, conditions $C_j'$, which characterize many prime numbers simultaneously. Afterwards, from the conditions $c_i, c_j'$, using the General Theorem again, we find new conditions $c_h''$ which characterize prime numbers simultaneously. And this method can be continued analogically.

<u>Remarks.</u>

Let $m_i = 1$ and $c_i$ represent the Simionov's theorem for all i.

(a) If D = 1 it results in V. Popa's theorem, which generalizes in its turn Cucurezeanu's theorem and the last one generalizes in its turn Clement's theorem!

(b) If D = $P/p_2$ and choosing conveniently the parameters $a_i$, $k_i$ for i = 1, 2, 3, it results in S. Patrizio's theorem.

**Several EXAMPLES:**

1. Let $p_1, p_2, \ldots, p_n$ be positive integers > 1, coprime integers two by two, and $1 \leq k_i \leq p_i$ for all i. Then:

$p_1, p_2, \ldots, p_n$ are simultaneously prime if and only if:

(T) $$\sum_{i=1}^{n} [(p_i-k_i)!(k_i-1)! - (-1)^{k_i}] \cdot \prod_{j \neq i} p_j \equiv 0 \pmod{p_1 p_2 \ldots p_n}$$

or

(U) $$\left(\sum_{i=1}^{n} [(p_i-k_i)!(k_i-1)!-(-1)^{k_i}] \cdot \prod_{j \neq i} p_j\right) / (p_{s+1} \ldots p_n) \equiv 0 \pmod{p_1 \ldots p_s}$$

or

(V) $$\sum_{i=1}^{n} [(p_i-k_i)!(k_i-1)!-(-1)^{k_i}] p_j/p_i \equiv 0 \pmod{p_j}$$

or

(W) $\quad \sum_{i=1}^{n} [(p_i-k_i)!(k_i-1)!-(-1)^{k_i}]/p_i$ is an integer.

2. Another relation example (using the first theorem from [4]): p is a prime positive integer if and only if $(p-3)!-(p-1)/2 \equiv 0 \pmod{p}$.

$$\sum_{i=1}^{n} [(p_i-3)!-(p_i-1)/2]p_i/p_i \equiv 0 \pmod{p_i}.$$

3. The odd numbers p and p + 2 are twin prime if and only if:

$(p-1)!(3p+2) + 2p + 2 \equiv \pmod{p(p+2)}$

or

$(p-1)!(p-2)-2 \equiv 0 \pmod{p(p+2)}$

or

$[(p-1)! + 1] / p + [(p-1)! \, 2 + 1] / (p+2)$ is an integer.

These twin prime characterizations differ from Clement's theorem $((p-1)!4 + p + 4 \equiv 0 \pmod{p(p+2)})$.

4. Let $(p,p+k) \sim 1$, then: p and p + k are prime simultaneously if and only if $(p-1)!(p+k) + (p+k-1)! \, p +$

$2p + k \equiv 0 \pmod{p(p+k)}$, which differs from I. Cucurezeanu's theorem ([1], p. 165): $k \cdot k! [(p-1)!+1] + [k!-(-1)^k] p \equiv 0 \pmod{p(p+k)}$).

5. Look at a characterization of a quadruple of primes for $p, p + 2, p + 6, p + 8$: $[(p-1)!+1]/p + [(p-1)!2!+1]/(p+2) + [(p-1)!6!+1]/(p+6) + [(p-1)!8!+1]/(p+8)$ be an integer.

6. For $p - 2, p, p + 4$ coprime integers two by two, we find the relation: $(p-1)!+p[(p-3)!+1]/(p-2)+p[(p+3)!+1]/(p+4) \equiv -1 \pmod{p}$, which differ from S. Patrizio's theorem $(8[(p+3)!/(p+4)] + 4[(p-3)!/(p-2)] \equiv -11 \pmod{p})$.

<u>References:</u>

[1]  Cucurezeanu, I., "Probleme de aritmetica si teoria numerelor", Ed. Tehnica, Bucuresti, 1966.

[2]  Patrizio, Serafino, "Generalizzazione del teorema di Wilson alle terne prime", Enseignement Math., Vol. 22(2), nr. 3-4, pp. 175-184, 1976.

[3]  Popa, Valeriu, "Asupra unor generalizari ale teoremei lui Clement", Studii si cercetari matematice, Vol. 24, Nr.

9, pp. 1435-1440, 1972.

[4]   Smarandache, Florentin, "Criterii ca un numar natural sa fie prim", Gazeta Matematica, Nr. 2, pp. 49-52; 1981; see Mathematical Reviews (USA): 83a: 10007.